\documentclass[11pt,reqno]{amsart}
\usepackage{amsfonts,amssymb,amsmath,color}

\DeclareMathOperator{\bmin}{\mathbf{min}}

\newtheorem{theorem}{Theorem}

\newtheorem{corollary}[theorem]{Corollary}
\newtheorem{lemma}[theorem]{Lemma}
\newtheorem{proposition}[theorem]{Proposition}

\newtheorem{remark}[theorem]{Remark}
\newtheorem{example}[theorem]{Example}

\numberwithin{equation}{section} \numberwithin{theorem}{section}

\begin{document}

\title{Pattern Recognition on Oriented Matroids: Layers of Tope Committees}

\author{Andrey O. Matveev}
\address{Data-Center Co., RU-620034, Ekaterinburg,
P.O.~Box~5, Russian~Federation} \email{aomatveev@\{dc.ru,
hotmail.com\}}


\begin{abstract}
A tope committee $\mathcal{K}^{\ast}$ for a simple oriented matroid $\mathcal{M}$ is a~subset of its maximal
covectors such that every positive halfspace of~$\mathcal{M}$ contains more than half of the covectors from
$\mathcal{K}^{\ast}$. The~structures of the family of all committees for $\mathcal{M}$, and of the family of
its committees that contain no pairs of opposites, are described. A~Farey subsequence associated with the
elements of the $m$th layer of the Boolean lattice of rank~$2m$ is explored.
\end{abstract}

\maketitle

\pagestyle{myheadings}

\markboth{PATTERN RECOGNITION ON ORIENTED MATROIDS}{A.O.~MATVEEV}

\thispagestyle{empty}

\tableofcontents

\section{Introduction}

Let $\boldsymbol{\mathcal{H}}:=\{\pmb{H}_1,\ldots,\pmb{H}_t\}$ be a finite {\em arrangement of\/} oriented
affine {\em hyperplanes} $\pmb{H}_i:=\{\pmb{x}\in\mathbb{R}^n:\ \langle\pmb{e}_i,\pmb{x}\rangle=a_i\}$ in
$\mathbb{R}^n$, where any two distinct vectors $\pmb{e}_i$ and $\pmb{e}_j$ are linearly independent,
$a_i\in\mathbb{R}$, and $\langle\pmb{e}_i,\pmb{x}\rangle:=\sum_{j=1}^n e_{ij}x_j$. See, e.g., \cite{OT,
St-HA} on hyperplane arrangements. The {\em regions\/} (or {\em chambers}) of $\boldsymbol{\mathcal{H}}$ are
the connected components of $\mathbb{R}^n-\bigcup_{\pmb{H}\in\boldsymbol{\mathcal{H}}}\pmb{H}$. A region
$\boldsymbol{R}$ lies on the {\em positive side\/} of the hyperplane $\pmb{H}_i$ if
$\langle\pmb{e}_i,\pmb{v}\rangle>a_i$, for a vector $\pmb{v}\in\pmb{R}$. Let $\boldsymbol{\mathcal{T}}_i^+$
denote the set of all regions of the arrangement $\boldsymbol{\mathcal{H}}$ that lie on the positive side of
the hyperplane $\pmb{H}_i$. A {\em majority committee of regions\/} (or a {\em committee}, for short) for the
arrangement $\boldsymbol{\mathcal{H}}$ is a subset of regions
$\boldsymbol{\mathcal{K}}^{\ast}:=\{\pmb{R}_1,\ldots,\pmb{R}_l\}$ such that it holds
$|\boldsymbol{\mathcal{K}}^{\ast}\cap\boldsymbol{\mathcal{T}}_i^+|>
\tfrac{|\boldsymbol{\mathcal{K}}^{\ast}|}{2}$, $1\leq i\leq t$. A representative system
$\{\pmb{w}_k\in\pmb{R}_k:\ 1\leq k\leq l\}$ is called a {\em committee\/} for the system of strict linear
inequalities $\{\langle\pmb{e}_i,\mathbf{x}\rangle>a_i:\ 1\leq i\leq t\}$, see,
e.g.,~\cite{AK-IEEE,AK-AMS,KhMR,M,MKh2}.

Oriented matroids are defined by various equivalent axiom systems,
and they can be thought of as a combinatorial abstraction of point
configurations over the reals, of real hyperplane arrangements, of
convex polytopes, and of directed graphs. Oriented matroids are
thoroughly studied, e.g, in~\cite{BK,BLSWZ,JGB,R-GZ,S,Z}.

An application of constructions, that generalize committees for
arrangements of linear hyperplanes, to the pattern recognition
problem in its abstract setting is as follows.

For\hfill a\hfill positive\hfill integer\hfill $t$,\hfill let\hfill $E_t$\hfill and\hfill $[1,t]$\hfill denote\hfill the\hfill set\hfill
$\{1,2,...,t\}$.\\ Let $\mathcal{M}:=(E_t,\mathcal{T})$ be a {\em
simple oriented matroid\/} (throughout the paper, `simple' means,
in a partly nonstandard way, that $\mathcal{M}$ has no
loops, parallel or {\sl antiparallel\/} elements) of rank
$r(\mathcal{M})\geq 2$, on the {\em ground set\/} $E_t$, with set
of {\em topes\/} $\mathcal{T}\subset\{-,+\}^{E_t}$.

The number of topes $|\mathcal{T}|$ can be computed with the help of the technique which\hfill was\hfill independently\hfill
proposed\hfill by\hfill Las~Vergnas\hfill and\hfill Zaslavsky,\\ see~\cite[\S{}4.6]{BLSWZ}. The {\em positive halfspace\/} associated
with an element $e\in E_t$ is the set $\mathcal{T}_e^+:=\{T\in\mathcal{T}:\ T(e)=+\}$. A subset
$\mathcal{K}^{\ast}\subset\mathcal{T}$ is called a {\em tope committee\/} for $\mathcal{M}$ if it holds
$|\mathcal{K}^{\ast}\cap\mathcal{T}_e^+|>\tfrac{|\mathcal{K}^{\ast}|}{2}$, for all $e\in E_t$. This paper is
a sequel to~\cite{M-Existence}, where it was shown that the family of tope committees for~$\mathcal{M}$,
denoted by $\mathbf{K}^{\ast}(\mathcal{M})$, is nonempty.

Let $\widetilde{\mathcal{M}}$ denote the {\em nontrivial
extension\/} of $\mathcal{M}$ by a new element $g$ which is not a
loop, and which is parallel or antiparallel to neither of the
elements of $E_t$; let $\sigma$ be the corresponding {\em
localization}. Fix a tope committee $\mathcal{K}^{\ast}$
for~$\mathcal{M}$. Let $\mathcal{C}^{\ast}$ denote the set of
cocircuits of $\mathcal{M}$, and suppose that the sets
$\bigl\{\bigl(X,\Sigma_K:=\sigma(X)\bigr):\
\mbox{$X\in\mathcal{C}^{\ast}$},\text{\ $X$ restriction of
$K$}\bigr\}$ are conformal, for all topes
$K\in\mathcal{K}^{\ast}$. The {\em committee decision rule},
corresponding to $\mathcal{K}^{\ast}$, relates the element $g$ to
a {\em class\/} $\mathbf{A}$ if $|\{K\in\mathcal{K}^{\ast}:\
\Sigma_K=+\}|<\tfrac{|\mathcal{K}^{\ast}|}{2}$; on the contrary, $g$ is {\em
recognized\/} as an element of the other {\em class\/}
$\mathbf{B}$ if $|\{K\in\mathcal{K}^{\ast}:\
\Sigma_K=+\}|>\tfrac{|\mathcal{K}^{\ast}|}{2}$,
see~\cite{M-Existence} for more on the two-class problem of
pattern recognition on oriented matroids.

In this paper, we present the structural description of the family $\mathbf{K}^{\ast}(\mathcal{M})$ of all
tope committees for a simple oriented matroid $\mathcal{M}$ which is not acyclic. The description involves
certain specific subsequences of the Farey sequences. See, e.g.,~\cite[Chapter~4]{GKP} on the standard {\em
Farey sequences\/} $\mathcal{F}_n$ {\em of order $n$}, which are defined to be the ascending sequences of
irreducible fractions $\tfrac{h}{k}$ such that $\tfrac{0}{1}\leq\tfrac{h}{k}\leq\tfrac{1}{1}$ and $k\leq n$.

Among interesting subsequences of $\mathcal{F}_n$ there are some sequences which have a neat set-theoretic
and combinatorial meaning: Let $A$ be a proper $m$-subset of a nonempty finite set $C$ of cardinality $n$.
For all nonempty subsets $B\subseteq C$, arrange the fractions $\tfrac{|B\cap A|}{|B|}$, reduced to their
lowest terms, in ascending order, without repetition. The resulting sequence
$\mathcal{F}\bigl(\mathbb{B}(n),m\bigr)$ (considered in~\cite[\S{}7]{AM}, where it was denoted by
$\mathcal{F}\bigl(\mathbb{B}(n),m;\rho\bigr)$, with $\rho$ making reference to the poset rank function on the
Boolean lattice $\mathbb{B}(n)$ of rank~$n$) inherits many properties of $\mathcal{F}_n$.

Since the cardinality of every positive halfspace of a simple
oriented matroid $\mathcal{M}:=(E_t,\mathcal{T})$ equals
$\tfrac{|\mathcal{T}|}{2}$, the Farey subsequences of special
interest are those analogous to
$\mathcal{F}\bigl(\mathbb{B}(|\mathcal{T}|),\tfrac{|\mathcal{T}|}{2}\bigr)$.

Although, for any positive integer $m$, the sequence $\mathcal{F}\bigl(\mathbb{B}(2m),m\bigr)$ is endowed
with the order-reversing bijection $\tfrac{h}{k}\mapsto\tfrac{k-h}{k}$, the entries within the~left
halfsequence
$\mathcal{F}^{\leq\frac{1}{2}}\bigl(\mathbb{B}(2m),m\bigr):=\bigl(f\in\mathcal{F}\bigl(\mathbb{B}(2m),m\bigr):\
f\leq\tfrac{1}{2}\bigr)$ and the~right halfsequence
$\mathcal{F}^{\geq\frac{1}{2}}\bigl(\mathbb{B}(2m),m\bigr):=\bigl(f\in\mathcal{F}\bigl(\mathbb{B}(2m),m\bigr):\
f\geq\tfrac{1}{2}\bigr)$ of $\mathcal{F}\bigl(\mathbb{B}(2m),m\bigr)$ exhibit different behavior. For
example, if $m>1$ then the numerators $h$ of the fractions $\tfrac{h}{k}\geq\tfrac{1}{2}$ are symmetrically
distributed with respect to the numerator~$2$ of the fraction $\tfrac{2}{3}$ which occupies the central
position in the subsequence $\mathcal{F}^{\geq\frac{1}{2}}\bigl(\mathbb{B}(2m),m\bigr)$; the subsequence
$\mathcal{F}^{\leq\frac{1}{2}}\bigl(\mathbb{B}(2m),m\bigr)$ does not have such a property. It is shown
in~\cite[\S{}3, \S{}4]{M-Integers} relying on results from the present paper that in fact the sequence of
numerators from $\mathcal{F}^{\leq\frac{1}{2}}\bigl(\mathbb{B}(2m),m\bigr)$ is exactly the sequence of
numerators from $\mathcal{F}_m$, and the sequence of numerators from
$\mathcal{F}^{\geq\frac{1}{2}}\bigl(\mathbb{B}(2m),m\bigr)$ is the sequence of denominators from
$\mathcal{F}_m$:

{\em Let $m>1$. The maps
\begin{align*}
\mathcal{F}^{\leq\frac{1}{2}}\!\bigl(\mathbb{B}(2m),m\bigr)&\to\mathcal{F}_m\ , &
\tfrac{h}{k}&\mapsto\tfrac{h}{k-h}\ ,\\
\mathcal{F}_m&\to\mathcal{F}^{\leq\frac{1}{2}}\!\bigl(\mathbb{B}(2m),m\bigr)\ , &
\tfrac{h}{k}&\mapsto\tfrac{h}{k+h}\ ,\\
\mathcal{F}^{\geq\frac{1}{2}}\!\bigl(\mathbb{B}(2m),m\bigr)&\to\mathcal{F}_m\ , &
\tfrac{h}{k}&\mapsto\tfrac{2h-k}{h}\ ,\intertext{and}
\mathcal{F}_m&\to\mathcal{F}^{\geq\frac{1}{2}}\!\bigl(\mathbb{B}(2m),m\bigr)\ , &
\tfrac{h}{k}&\mapsto\tfrac{k}{2k-h}\ ,
\end{align*}
are order-preserving and bijective.

The maps
\begin{align*}
\mathcal{F}^{\leq\frac{1}{2}}\!\bigl(\mathbb{B}(2m),m\bigr)&\to\mathcal{F}_m\ , &
\tfrac{h}{k}&\mapsto\tfrac{k-2h}{k-h}\ ,\\
\mathcal{F}_m&\to\mathcal{F}^{\leq\frac{1}{2}}\!\bigl(\mathbb{B}(2m),m\bigr)\ , &
\tfrac{h}{k}&\mapsto\tfrac{k-h}{2k-h}\ ,\\
\mathcal{F}^{\geq\frac{1}{2}}\!\bigl(\mathbb{B}(2m),m\bigr)&\to\mathcal{F}_m\ , &
\tfrac{h}{k}&\mapsto\tfrac{k-h}{h}\ ,\intertext{and}
\mathcal{F}_m&\to\mathcal{F}^{\geq\frac{1}{2}}\!\bigl(\mathbb{B}(2m),m\bigr)\ , &
\tfrac{h}{k}&\mapsto\tfrac{k}{k+h}\ ,
\end{align*}
are order-reversing and bijective.
}

If $\mathcal{K}^{\ast}$ is a tope committee for an oriented matroid $\mathcal{M}:=(E_t,\mathcal{T})$ then the
disjoint union $\mathcal{K}^{\ast}\dot\cup\{T,-T\}$ of $\mathcal{K}^{\ast}$ with a pair of opposites
$\{T,-T\}\subset\mathcal{T}$ is also a committee for $\mathcal{M}$. In a similar way, if
$\mathcal{K}'^{\ast}$ and $\mathcal{K}''^{\ast}$ are disjoint committees, then their union is a committee as
well. Such redundant committees are not of applied importance because in practice one searches in reverse
direction, for inclusion-minimal committees: if $\mathcal{K}'^{\ast}$ and $\mathcal{K}''^{\ast}$ are tope
committees for $\mathcal{M}$, with $\mathcal{K}'^{\ast}\subsetneqq\mathcal{K}''^{\ast}$, then the committee
$\mathcal{K}'^{\ast}$ is preferred. A committee $\mathcal{K}^{\ast}$ is called {\em minimal\/} if any its
proper subset is not a committee. Minimal committees do not contain opposites. Among the minimal committees
the most interesting are the committees of minimal cardinality, so-called {\em minimum committees}.

In Section~\ref{Relative-Blocking} of the paper, we interpret tope
committees for oriented matroids in terms of relatively blocking
elements in Boolean lattices. The~Farey subsequences
$\mathcal{F}\bigl(\mathbb{B}(2m),m\bigr)$ are explored in
Section~\ref{Section-Farey-2m-m}. Subfamilies of the~family
$\mathbf{K}^{\ast}(\mathcal{M})$ of tope committees for an
oriented matroid $\mathcal{M}:=(E_t,\mathcal{T})$, lying on layers
of the Boolean lattice $\mathbb{B}(\mathcal{T})$ of subsets of the
tope set $\mathcal{T}$, are described in
Section~\ref{section:layers}. In
Section~\ref{Section-No-Opposites}, we consider the subfamily
$\overset{\circ}{\mathbf{K}}{}^{\ast}(\mathcal{M})\subset
\mathbf{K}^{\ast}(\mathcal{M})$ of tope committees for
$\mathcal{M}$, that contain no pairs of opposites, and we describe
its structure.

\section{Tope Committees and Relative Blocking}
\label{Relative-Blocking}

Let $\mathcal{M}=(E_t,\mathcal{T})$ be a simple oriented matroid on the ground set $E_t$, with set of topes
$\mathcal{T}$. The Boolean lattice of all subsets of $\mathcal{T}$ is denoted by~$\mathbb{B}(\mathcal{T})$.
Let $\Upsilon$ denote the antichain $\{\upsilon_1,\ldots,\upsilon_t\}$ corresponding to the~family of
positive halfspaces $\{\mathcal{T}_1^+,\ldots,\mathcal{T}_t^+\}$ which is thought of as a subset of~$\mathbb{B}(\mathcal{T})$. The~poset rank $\rho(\upsilon_e)$ of any element $\upsilon_e$ in
$\mathbb{B}(\mathcal{T})$ is $\tfrac{|\mathcal{T}|}{2}$, therefore the antichain $\Upsilon$ is {\em pure} in
the sense that it lies entirely on a layer of the graded lattice $\mathbb{B}(\mathcal{T})$.

Interpret a tope committee for $\mathcal{M}$ as an element
$b\in\mathbb{B}(\mathcal{T})$. Then the~family
\begin{equation*}
\mathbf{K}^{\ast}(\mathcal{M}):=\left\{\mathcal{K}^{\ast}\subset\mathcal{T}:
\ |\mathcal{K}^{\ast}\cap\mathcal{T}_e^+|>
\tfrac{1}{2}|\mathcal{K}^{\ast}|,\ \forall e\in E_t\right\}
\end{equation*}
of all tope committees for $\mathcal{M}$ can be viewed as the
subposet
\begin{equation*}
\mathbf{I}_{\frac{1}{2}}\bigl(\mathbb{B}(\mathcal{T}),
\Upsilon\bigr):=\left\{b\in\mathbb{B}(\mathcal{T}):\ \rho(b)>0,\
\tfrac{\rho(b\wedge\upsilon_e)}{\rho(b)}>\tfrac{1}{2},\ \forall
e\in E_t \right\}
\end{equation*}
of all {\em relatively $\tfrac{1}{2}$-blocking elements\/} for the
antichain $\Upsilon$ in $\mathbb{B}(\mathcal{T})$. Relative
blocking in posets is discussed in~\cite{AM}. The antichain
$\mathfrak{y}_{\tfrac{1}{2}}\bigl(\mathbb{B}(\mathcal{T}),
\Upsilon\bigr):=
\bmin\mathbf{I}_{\frac{1}{2}}\bigl(\mathbb{B}(\mathcal{T}),
\Upsilon\bigr)$, called in~\cite{AM} the {\em relative
$\tfrac{1}{2}$-blocker of\/ $\Upsilon$ in
$\mathbb{B}(\mathcal{T})$}, is the~family of all minimal tope
committees for $\mathcal{M}$; throughout the paper,~$\bmin$ denotes the set of all minimal elements of a subposet.

For a graded poset $\mathfrak{P}$ with rank function $\rho$, and
for a nonnegative integer~$k$, we let $\mathfrak{P}^{(k)}$ denote
the $k$th {\em layer\/} of $\mathfrak{P}$, that is the antichain
$\{p\in\mathfrak{P}:\ \rho(p)=k\}$. If $A$ is an antichain in
$\mathfrak{P}$, then $\mathfrak{I}(A)$ and $\mathfrak{F}(A)$
denote the order ideal and filter in $\mathfrak{P}$, generated by
$A$, respectively.

For $k\in[1,|\mathcal{T}|-1]$, the subposet
\begin{equation}
\label{eq:10}
\mathbf{I}_{\frac{1}{2},k}\bigl(\mathbb{B}(\mathcal{T}),
\Upsilon\bigr):=\mathbb{B}(\mathcal{T})^{(k)}\ \cap\
\mathbf{I}_{\frac{1}{2}}\bigl(\mathbb{B}(\mathcal{T}),
\Upsilon\bigr)
\end{equation}
of elements of rank $k$ from
$\mathbf{I}_{\frac{1}{2}}\bigl(\mathbb{B}(\mathcal{T}),
\Upsilon\bigr)$ is the antichain
\begin{equation*}
\mathbb{B}(\mathcal{T})^{(k)}\ \cap\ \bigcap_{e\in
E_t}\mathfrak{F}\!\left(\mathfrak{I}(\upsilon_e)\ \cap\
\mathbb{B}(\mathcal{T})^{(\lceil(k+1)/2\rceil)}\right)\ ,
\end{equation*}
see~\cite[Proposition~5.1(ii)]{AM}.

\section{The Farey Subsequence $\mathcal{F}\bigl(\mathbb{B}(2m),m\bigr)$}
\label{Section-Farey-2m-m}

Let $C$ be a finite nonempty set of even cardinality $2m$, and $A$ an $m$-subset of $C$. Arrange in ascending
order (without repetition) the fractions $\tfrac{|B\cap A|}{|B|}$, reduced to their lowest terms, for all
nonempty subsets $B\subseteq C$; the resulting Farey subsequence is
\begin{equation*}
\mathcal{F}\bigl(\mathbb{B}(2m),m\bigr):=
\bigl(\tfrac{h}{k}\in\mathcal{F}_{2m}:\ h\leq m,\ k-h\leq m
\bigr)\ .
\end{equation*}

\begin{example}
\label{prop:6}
\begin{equation*}
\mathcal{F}\bigl(\mathbb{B}(8),4\bigr)=\Bigl(\tfrac{0}{1}<\tfrac{1}{5}
<\tfrac{1}{4}<\tfrac{1}{3}<\tfrac{2}{5}<
\tfrac{3}{7}<\tfrac{1}{2}<\tfrac{4}{7}<\tfrac{3}{5}<\tfrac{2}{3}<
\tfrac{3}{4}<\tfrac{4}{5}<\tfrac{1}{1}\Bigr)\ ;
\end{equation*}
\begin{multline*}
\mathcal{F}\bigl(\mathbb{B}(10),5\bigr)=\Bigl(\tfrac{0}{1}<
\tfrac{1}{6}<\tfrac{1}{5}<\tfrac{1}{4}<\tfrac{2}{7}<\tfrac{1}{3}<\tfrac{3}{8}<
\tfrac{2}{5}<\tfrac{3}{7}<\tfrac{4}{9}\\
<\tfrac{1}{2}<\tfrac{5}{9}<\tfrac{4}{7}<
\tfrac{3}{5}<\tfrac{5}{8}<\tfrac{2}{3}<\tfrac{5}{7}<\tfrac{3}{4}<\tfrac{4}{5}<
\tfrac{5}{6}<\tfrac{1}{1}\Bigr)\ .
\end{multline*}
\end{example}

In this section, we explore such sequences, and we start by recalling some basic properties of general
sequences
\begin{equation*}
\mathcal{F}\bigl(\mathbb{B}(n),m\bigr):=\bigl(\tfrac{h}{k}\in\mathcal{F}_{n}:\ h\leq m,\ k-h\leq n-m\bigr)\ ,
\end{equation*}
see also~\cite{{M-Integers-II}}.

\begin{lemma}\cite[Proposition~7.5]{AM}
\label{prop:1} Let
$\tfrac{h}{k}\in\mathcal{F}\bigl(\mathbb{B}(n),m\bigr)
-\{\tfrac{0}{1},\tfrac{1}{1}\}$,
where $0<m<n$.
\begin{itemize}
\item[\rm(i)]
Let $x_0$ be the integer such that $kx_0\equiv -1\pmod{h}$ and
$m-h+1\leq x_0\leq m$. Define integers $y_0$ and $t^{\ast}$ by
$y_0:=\tfrac{kx_0+1}{h}$ and
$t^\ast:=\left\lfloor\min\{\tfrac{m-x_0}{h},\tfrac{n-y_0}{k},
\tfrac{n-m+x_0-y_0}{k-h}\}\right\rfloor$.

The fraction $\tfrac{x_0+t^\ast h}{y_0+t^\ast k}$ precedes the
fraction $\tfrac{h}{k}$ in
$\mathcal{F}\bigl(\mathbb{B}(n),m\bigr)$.

\item[\rm(ii)]
Let $x_0$ be the integer such that $kx_0\equiv 1\pmod{h}$ and
$m-h+1\leq x_0\leq m$. Define integers $y_0$ and $t^{\ast}$ by
$y_0:=\tfrac{kx_0-1}{h}$ and
$t^{\ast}:=\left\lfloor\min\{\tfrac{m-x_0}{h},\tfrac{n-y_0}{k},
\tfrac{n-m+x_0-y_0}{k-h}\}\right\rfloor$.

The fraction $\tfrac{x_0+t^{\ast} h}{y_0+t^{\ast} k}$ succeeds the
fraction $\tfrac{h}{k}$ in
$\mathcal{F}\bigl(\mathbb{B}(n),m\bigr)$.
\end{itemize}
\end{lemma}

\newpage

\begin{lemma}\cite[Proposition~7.8]{AM}
\label{prop:2} Let
$\frac{h_j}{k_j}<\frac{h_{j+1}}{k_{j+1}}<\frac{h_{j+2}}{k_{j+2}}$
be three successive fractions of
$\mathcal{F}\bigl(\mathbb{B}(n),m\bigr)$, where $0<m<n$.

\begin{itemize}
\item[\rm(i)] The integers $h_j$ and $k_j$ are computed by
{\small
\begin{align*}
h_j&=\left\lfloor\min\left\{\frac{h_{j+2}+m}{h_{j+1}},
\frac{k_{j+2}+n}{k_{j+1}},\frac{k_{j+2}
-h_{j+2}+n-m}{k_{j+1}-h_{j+1}}\right\}\right\rfloor
h_{j+1}-h_{j+2}\ ,\\
k_j&=\left\lfloor\min\left\{\frac{h_{j+2}+m}{h_{j+1}},
\frac{k_{j+2}+n}{k_{j+1}},\frac{k_{j+2}
-h_{j+2}+n-m}{k_{j+1}-h_{j+1}}\right\}\right\rfloor
k_{j+1}-k_{j+2}\ .
\end{align*}
}

\item[\rm(ii)] The integers $h_{j+2}$ and $k_{j+2}$ are computed by
{\small
\begin{align*}
h_{j+2}&=\left\lfloor\min\left\{\frac{h_j+m}{h_{j+1}},
\frac{k_j+n}{k_{j+1}},\frac{k_j
-h_j+n-m}{k_{j+1}-h_{j+1}}\right\}\right\rfloor h_{j+1}-h_j\ ,\\
k_{j+2}&=\left\lfloor\min\left\{\frac{h_j+m}{h_{j+1}},
\frac{k_j+n}{k_{j+1}},\frac{k_j
-h_j+n-m}{k_{j+1}-h_{j+1}}\right\}\right\rfloor k_{j+1}-k_j\ .
\end{align*}
}
\end{itemize}
\end{lemma}

The first observation is as follows:

\begin{remark}
\label{prop:4} The map
\begin{equation*}
\mathcal{F}\bigl(\mathbb{B}(2m),m\bigr)\to\mathcal{F}\bigl(\mathbb{B}(2m),m\bigr)\ ,\ \ \
\tfrac{h}{k}\mapsto\tfrac{k-h}{k}
\end{equation*}
is order-reversing and bijective.
\end{remark}

In the case where $n:=2m$, Lemmas~\ref{prop:1} and~\ref{prop:2} can be refined:

\begin{corollary}
\label{prop:3}
\begin{itemize}
\item[\rm(i)] Let
$\tfrac{h}{k}\in\mathcal{F}\bigl(\mathbb{B}(2m),m\bigr)$. Suppose
that $\tfrac{h}{k}>\tfrac{1}{2}$.

Let $x_0$ be the integer such that $kx_0\equiv -1\pmod{h}$ and
$m-h+1\leq x_0\leq m$. The fraction
\begin{equation*}
x_0\!\Big/ \tfrac{kx_0+1}{h}
\end{equation*}
precedes the fraction $\tfrac{h}{k}$ in
$\mathcal{F}\bigl(\mathbb{B}(2m),m\bigr)$.

\item[\rm(ii)] Let
$\tfrac{h}{k}\in\mathcal{F}\bigl(\mathbb{B}(2m),m\bigr)$. Suppose
that $\tfrac{1}{2}\leq\tfrac{h}{k}<\tfrac{1}{1}$.

Let $x_0$ be the integer such that $kx_0\equiv 1\pmod{h}$ and
$m-h+1\leq x_0\leq m$. The fraction
\begin{equation}
\label{eq:1} x_0\!\Big/ \tfrac{kx_0-1}{h}
\end{equation}
succeeds the fraction $\tfrac{h}{k}$ in
$\mathcal{F}\bigl(\mathbb{B}(2m),m\bigr)$. In particular, the
fraction $\tfrac{m}{2m-1}$ succeeds $\tfrac{1}{2}$.
\item[\rm(iii)] Let
$\frac{h_j}{k_j}<\frac{h_{j+1}}{k_{j+1}}<\frac{h_{j+2}}{k_{j+2}}$
be three successive fractions in
$\mathcal{F}\bigl(\mathbb{B}(2m),m\bigr)$, where $m>1$, with
$\tfrac{h_j}{k_j}\geq\tfrac{1}{2}$.

The integers $h_j$, $k_j$, $h_{j+2}$ and $k_{j+2}$ are computed by
{\small
\begin{align}
\label{eq:5}
h_j&=\left\lfloor\frac{h_{j+2}+m}{h_{j+1}}\right\rfloor
h_{j+1}-h_{j+2}\ , &
k_j&=\left\lfloor\frac{h_{j+2}+m}{h_{j+1}}\right\rfloor
k_{j+1}-k_{j+2}\ ,\\ \label{eq:6}
h_{j+2}&=\left\lfloor\frac{h_j+m}{h_{j+1}}\right\rfloor
h_{j+1}-h_j\ , &
k_{j+2}&=\left\lfloor\frac{h_j+m}{h_{j+1}}\right\rfloor
k_{j+1}-k_j\ .
\end{align}
}
\end{itemize}
\end{corollary}

\begin{proof} We prove~(ii); assertion~(i) is proved in a
similar way.

{\rm(ii)} In terms of Lemma~\ref{prop:1}(ii),
$t^{\ast}:=\left\lfloor\min\{\tfrac{m-x_0}{h},\tfrac{2m-y_0}{k},
\tfrac{m+x_0-y_0}{k-h}\}\right\rfloor$. Since $\tfrac{h}{k}\geq\tfrac{1}{2}$, we have
$\min\{\tfrac{m-x_0}{h},\tfrac{2m-y_0}{k},
\tfrac{m+x_0-y_0}{k-h}\}=\tfrac{m-x_0}{h}$. Therefore
$t^{\ast}=\left\lfloor\tfrac{m-x_0}{h}\right\rfloor$. But the
constraint $0\leq m-x_0\leq h-1$ implies $t^{\ast}=0$, and the
assertion follows from Lemma~\ref{prop:1}(ii).

{\rm(iii)} We prove~(\ref{eq:6}). Lemma~\ref{prop:2}(ii) implies
that {\small
\begin{align*}
h_{j+2}&=\left\lfloor\min\left\{\frac{h_j+m}{h_{j+1}},
\frac{k_j+2m}{k_{j+1}},\frac{k_j
-h_j+m}{k_{j+1}-h_{j+1}}\right\}\right\rfloor h_{j+1}-h_j\ ,\\
k_{j+2}&=\left\lfloor\min\left\{\frac{h_j+m}{h_{j+1}},
\frac{k_j+2m}{k_{j+1}},\frac{k_j
-h_j+m}{k_{j+1}-h_{j+1}}\right\}\right\rfloor k_{j+1}-k_j\ .
\end{align*}
} Since $\tfrac{h_j}{k_j}\geq\tfrac{1}{2}$, we have
$\min\left\{\frac{h_j+m}{h_{j+1}},
\frac{k_j+2m}{k_{j+1}},\frac{k_j
-h_j+m}{k_{j+1}-h_{j+1}}\right\}=\frac{h_j+m}{h_{j+1}}$;
hence~(\ref{eq:6}).
\end{proof}

Analogous properties of the left halfsequence of $\mathcal{F}\bigl(\mathbb{B}(2m),m\bigr)$ are presented in
Proposition~\ref{prop:7} below.

For a fraction $f:=\tfrac{h}{k}$, we denote by $\underline{f}:=h$
and $\overline{f}:=k$ its numerator and denominator, respectively.

The fractions $f_s\in\mathcal{F}\bigl(\mathbb{B}(n),m\bigr)$ are
always indexed starting with zero, thus~$f_0:=\tfrac{0}{1}$.

\begin{proposition}
\label{prop:5} Let\hfill
$f_s,f_t\in\mathcal{F}\bigl(\mathbb{B}(2m),m\bigr)$,\hfill where\hfill $m>1$,\hfill
$f_s:=\tfrac{1}{2}$\\ and $f_t:=\tfrac{2}{3}$. Fractions from
$\mathcal{F}\bigl(\mathbb{B}(2m),m\bigr)$ satisfy the equalities
\begin{align}
\label{eq:2} \underline{f_{t+v}}&=\underline{f_{t-v}}\ , \\
\label{eq:3}
\bigl(\,\overline{f_{t+v}}+\overline{f_{t-v}}\,\,\bigr) \Big/
\underline{f_{t+v}}\ &=\ 3\ .
\end{align}
for all $v$, $0\leq v\leq t-s$.
\end{proposition}

\begin{proof}
Notice that~(\ref{eq:3}) holds for $v:=0$.

Let $v:=1$. By Corollary~\ref{prop:3}(ii), the fraction
$\tfrac{h_t}{k_t}:=f_t$ precedes the fraction
$\tfrac{h_{t+v}}{k_{t+v}}:=f_{t+v}$ with $h_{t+v}=x_0$ and
$k_{t+v}=\tfrac{k_tx_0-1}{h_t}$, where $3x_0\equiv 1\pmod{2}$ and
$m-1\leq x_0\leq m$; therefore the numerator and denominator of
$f_{t+v}$ are
\begin{align*}
h_{t+v}:=\underline{f_{t+v}}=x_0&=\begin{cases}m-1,&\text{if $m$
is even},\\ m,&\text{if $m$ is odd};\end{cases}\\
k_{t+v}:=\overline{f_{t+v}}=\frac{k_t
x_0-1}{h_t}&=\begin{cases}\frac{3m-4}{2},&\text{if $m$ is even},\\
\frac{3m-1}{2},&\text{if $m$ is odd}.\end{cases}
\end{align*}

Corollary~\ref{prop:3}(i) implies that the numerator and
denominator of the fraction $\tfrac{h_{t-v}}{k_{t-v}}:=f_{t-v}$
are
\begin{align*}
h_{t-v}:=\underline{f_{t-v}}&=\begin{cases}m-1,&\text{if $m$ is
even},\\ m,&\text{if $m$ is odd};\end{cases}\\
k_{t-v}:=\overline{f_{t-v}}&=\begin{cases}\frac{3m-2}{2},&\text{if
$m$ is even},\\ \frac{3m+1}{2},&\text{if $m$ is odd}.\end{cases}
\end{align*}

We see that equalities~(\ref{eq:2}) and~(\ref{eq:3}) hold for
$v:=1$.

Let $v:=2$. Equalities~(\ref{eq:6}) and~(\ref{eq:5}) yield
\begin{align*}
h_{t+v}:=\underline{f_{t+v}}&=
\left\lfloor\frac{h_{t+(v-2)}+m}{h_{t+(v-1)}}\right\rfloor
h_{t+(v-1)}-h_{t+(v-2)}\ ,\\ h_{t-v}:=\underline{f_{t-v}}&=
\left\lfloor\frac{h_{t-(v-2)}+m}{h_{t-(v-1)}}\right\rfloor
h_{t-(v-1)}-h_{t-(v-2)}\ ,
\end{align*}
respectively. It has been shown that $h_{t+(v-1)}=h_{t-(v-1)}$
and, by convention, we have $h_{t+(v-2)}=h_{t-(v-2)}=2$.
Thus,~(\ref{eq:2}) holds for $v:=2$.

Let $k_{t+v}:=\overline{f_{t+v}}$ and
$k_{t-v}:=\overline{f_{t-v}}$. Equalities~(\ref{eq:6})
and~(\ref{eq:5}) yield
\begin{multline*}
\bigl(\,\overline{f_{t+v}}+\overline{f_{t-v}}\,\,\bigr) \Big/
\underline{f_{t+v}}=\frac{k_{t+v}+k_{t-v}}{h_{t+v}}\\=
\frac{\left\lfloor\frac{h_{t+(v-2)}+m}{h_{t+(v-1)}}\right\rfloor
k_{t+(v-1)}-k_{t+(v-2)} +
\left\lfloor\frac{h_{t-(v-2)}+m}{h_{t-(v-1)}}\right\rfloor
k_{t-(v-1)}-k_{t-(v-2)} }
{\left\lfloor\frac{h_{t+(v-2)}+m}{h_{t+(v-1)}}\right\rfloor
h_{t+(v-1)}-h_{t+(v-2)}}\\=
\frac{\left\lfloor\frac{h_{t+(v-2)}+m}{h_{t+(v-1)}}\right\rfloor
\left(k_{t+(v-1)}+k_{t-(v-1)}\right)-\left(k_{t+(v-2)}+k_{t-(v-2)}\right)
}{\left\lfloor\frac{h_{t+(v-2)}+m}{h_{t+(v-1)}}\right\rfloor
h_{t+(v-1)}-h_{t+(v-2)}}\\= \frac{
3\!\!\left\lfloor\frac{h_{t+(v-2)}+m}{h_{t+(v-1)}}\right\rfloor
h_{t+(v-1)}-3h_{t+(v-2)}}{\left\lfloor\frac{h_{t+(v-2)}+m}{h_{t+(v-1)}}\right\rfloor
h_{t+(v-1)}-h_{t+(v-2)}}=3\ ,
\end{multline*}
that is, we have obtained~(\ref{eq:3}) for $v:=2$.

Let $v$ be any integer $\geq 3$ such that the fraction $f_{t-v}$
from $\mathcal{F}\bigl(\mathbb{B}(2m),m\bigr)$ is greater than or
equal to $\tfrac{1}{2}$. Equalities~(\ref{eq:2}) and~(\ref{eq:3})
are proved by induction.

Recall that for $v:=t-s$, by convention, we have
$\underline{f_{t-v}}=1$ and $\overline{f_{t-v}}=2$. We conclude
from~(\ref{eq:2}) and~(\ref{eq:3}) that $f_{t+v}=\tfrac{1}{1}$.
\end{proof}

Along with Remark~\ref{prop:4}, Proposition~\ref{prop:5} leads to
the following observation:

\newpage

\begin{corollary}
\begin{itemize}
\item[\rm(i)]
The maps
\begin{align*}
\mathcal{F}^{\geq\frac{1}{2}}\bigl(\mathbb{B}(2m),m\bigr)&\to
\mathcal{F}^{\geq\frac{1}{2}}\bigl(\mathbb{B}(2m),m\bigr)\ , & \tfrac{h}{k}&\mapsto\tfrac{h}{3h-k}\
,\intertext{and} \mathcal{F}^{\leq\frac{1}{2}}\bigl(\mathbb{B}(2m),m\bigr)&\to
\mathcal{F}^{\leq\frac{1}{2}}\bigl(\mathbb{B}(2m),m\bigr)\ , & \tfrac{h}{k}&\mapsto\tfrac{k-2h}{2k-3h}\ ,
\end{align*}
are order-reversing and bijective.
\item[\rm(ii)] Let $f_s,f_t\in\mathcal{F}\bigl(\mathbb{B}(2m),m\bigr)$,
where $m>1$, $f_s:=\tfrac{1}{3}$ and $f_t:=\tfrac{1}{2}$.
Fractions from $\mathcal{F}\bigl(\mathbb{B}(2m),m\bigr)$ satisfy
the equalities
\begin{equation*}
\bigl(\,\overline{f_{s+v}}+\overline{f_{s-v}}\,\,\bigr) \Big/
\bigl(\,\underline{f_{s+v}}+\underline{f_{s-v}}\,\,\bigr)\ =\ 3\ ,
\end{equation*}
for all $v$, $0\leq v\leq t-s$.
\end{itemize}
\end{corollary}

Thus, fragments of $\mathcal{F}\bigl(\mathbb{B}(2m),m\bigr)$, for large $m$, look like
\begin{multline*}
\mathcal{F}\bigl(\mathbb{B}(2m),m\bigr)\ \approx\
\Bigl(\tfrac{0}{1}<\tfrac{1}{m+1}<\tfrac{1}{m}<\tfrac{1}{m-1}
<\cdots<\tfrac{1}{3}<\cdots
\\ <\tfrac{m-3}{2m-5}<\tfrac{m-2}{2m-3}<\tfrac{m-1}{2m-1}<\tfrac{1}{2}
<\tfrac{m}{2m-1}< \tfrac{m-1}{2m-3}< \tfrac{m-2}{2m-5}<\cdots\\
<\cdots<\tfrac{2}{3}<\cdots<\tfrac{m-2}{m-1}<\tfrac{m-1}{m}
<\tfrac{m}{m+1}<\tfrac{1}{1}\Bigr)\ ;
\end{multline*}
see Example~\ref{prop:6}; this observation is discussed in more detail in~\cite[\S{}4]{M-Neighbors}. In particular, the~fraction
$\tfrac{h}{k}$, with
\begin{equation*}
h=\begin{cases}\frac{m-2}{2},&\text{if $m$ is even},\\
\frac{m-1}{2},&\text{if $m$ is odd},\end{cases}\ \ \ \ \
k=\begin{cases}\frac{3m-4}{2},&\text{if $m$ is even},\\
\frac{3m-1}{2},&\text{if $m$ is odd},\end{cases}
\end{equation*}
precedes $\tfrac{1}{3}$ in
$\mathcal{F}\bigl(\mathbb{B}(2m),m\bigr)$, and the fraction
$\tfrac{h}{k}$, with
\begin{equation*}
h=\begin{cases}\frac{m}{2},&\text{if $m$ is even},\\
\frac{m+1}{2},&\text{if $m$ is odd},\end{cases}\ \ \ \ \
k=\begin{cases}\frac{3m-2}{2},&\text{if $m$ is even},\\
\frac{3m+1}{2},&\text{if $m$ is odd},\end{cases}
\end{equation*}
succeeds $\tfrac{1}{3}$ in
$\mathcal{F}\bigl(\mathbb{B}(2m),m\bigr)$; see~\cite{M-Neighbors} for more on neighboring fractions in~$\mathcal{F}\bigl(\mathbb{B}(2m),m\bigr)$.

We conclude this section by presenting an analogue of Corollary~\ref{prop:3}. It describes some properties
(see~\cite[\S{}4]{M-Integers}) of the left halfsequence of $\mathcal{F}\bigl(\mathbb{B}(2m),m\bigr)$:

\begin{proposition}
\label{prop:7}
\begin{itemize}
\item[\rm(i)]
Let $\tfrac{h}{k}\in\mathcal{F}\bigl(\mathbb{B}(2m),m\bigr)$. Suppose that
$\tfrac{0}{1}<\tfrac{h}{k}\leq\tfrac{1}{2}$. Let $x_0$ be the integer such that $hx_0\equiv 1\pmod{(k-h)}$
and $m-k+h+1\leq x_0\leq m$. The fraction
\begin{equation*}
\tfrac{hx_0-1}{k-h}\!\!\Bigm/\!\!\tfrac{kx_0-1}{k-h}
\end{equation*}
precedes $\tfrac{h}{k}$ in $\mathcal{F}\bigl(\mathbb{B}(2m),m\bigr)$. In particular, the fraction
$\tfrac{m-1}{2m-1}$ precedes $\tfrac{1}{2}$.

\item[\rm(ii)]
Let $\tfrac{h}{k}\in\mathcal{F}\bigl(\mathbb{B}(2m),m\bigr)$. Suppose that
$\tfrac{0}{1}\leq\tfrac{h}{k}<\tfrac{1}{2}$. Let $x_0$ be the integer such that $hx_0\equiv -1\pmod{(k-h)}$
and $m-k+h+1\leq x_0\leq m$. The fraction
\begin{equation*}
\tfrac{hx_0+1}{k-h}\!\!\Bigm/\!\!\tfrac{kx_0+1}{k-h}
\end{equation*}
succeeds $\tfrac{h}{k}$ in $\mathcal{F}\bigl(\mathbb{B}(2m),m\bigr)$.

\item[\rm(iii)] Let
$\frac{h_j}{k_j}<\frac{h_{j+1}}{k_{j+1}}<\frac{h_{j+2}}{k_{j+2}}$
be three successive fractions in
$\mathcal{F}\bigl(\mathbb{B}(2m),m\bigr)$, where $m>1$, with
$\tfrac{h_{j+2}}{k_{j+2}}\leq\tfrac{1}{2}$.

The integers $h_j$, $k_j$, $h_{j+2}$ and $k_{j+2}$ are computed by
{\small
\begin{align*}
h_j&=\left\lfloor\frac{k_{j+2}-h_{j+2}+m}{k_{j+1}-h_{j+1}}\right\rfloor
h_{j+1}-h_{j+2}\ , &
k_j&=\left\lfloor\frac{k_{j+2}-h_{j+2}+m}{k_{j+1}-h_{j+1}}\right\rfloor
k_{j+1}-k_{j+2}\ ,\\
h_{j+2}&=\left\lfloor\frac{k_j-h_j+m}{k_{j+1}-h_{j+1}}\right\rfloor
h_{j+1}-h_j\ , &
k_{j+2}&=\left\lfloor\frac{k_j-h_j+m}{k_{j+1}-h_{j+1}}\right\rfloor
k_{j+1}-k_j\ .
\end{align*}
}
\end{itemize}
\end{proposition}

\section{Layers of Tope Committees}

\label{section:layers}

We now describe the structure of the family $\mathbf{K}^{\ast}(\mathcal{M})$ of all tope committees for an
oriented matroid $\mathcal{M}:=(E_t,\mathcal{T})$ which is not acyclic.

The following assertion is a consequence
of~\cite[Theorem~8.4]{AM}.

\begin{proposition} \label{prop:8}
Let $\mathcal{M}$ be a simple oriented matroid, which is not
acyclic, on the ground set $E_t$, with set of topes $\mathcal{T}$.

On the one hand,
\begin{equation*}
\mathbf{I}_{\frac{1}{2}}\bigl(\mathbb{B}(\mathcal{T}),
\Upsilon\bigr)\\ = \dot{\bigcup_{3\leq k\leq|\mathcal{T}|-3}}\biggl(\
\mathbb{B}(\mathcal{T})^{(k)}\ \, \cap\ \bigcap_{e\in E_t
}\mathfrak{F}\!\left(\mathfrak{I}(\upsilon_e)\ \cap\
\mathbb{B}(\mathcal{T})^{(\lceil(k+1)/2\rceil)}\right)\ \biggr)
\end{equation*}
and, in particular,
\begin{equation*}
\mathbf{I}_{\frac{1}{2},3}\bigl(\mathbb{B}(\mathcal{T}),
\Upsilon\bigr) =
\mathbb{B}(\mathcal{T})^{(3)}\cap\bigcap_{e\in E_t
}\mathfrak{F}\!\left(\mathfrak{I}(\upsilon_e)\ \cap\
\mathbb{B}(\mathcal{T})^{(2)}\right)\ .
\end{equation*}
On the other hand,
\begin{multline*}
\mathbf{I}_{\frac{1}{2}}\bigl(\mathbb{B}(\mathcal{T}),
\Upsilon\bigr)\ \ \ =\ \ \ \bigcap_{e\in E_t}\ \ \
\bigcup_{f\in\mathcal{F}(\mathbb{B}(|\mathcal{T}|),\frac{|\mathcal{T}|}{2}):\
\frac{1}{2}<f}\ \ \ \
\bigcup_{s\in\left[1,\left\lfloor|\mathcal{T}|/(2\underline{f})\right\rfloor\right]}
\\ \Biggl(\
\mathbb{B}(\mathcal{T})^{(s\cdot\overline{f})}\ \ \cap \ \
\biggl(\ \mathfrak{F}\!\left(\mathfrak{I}(\upsilon_e)\ \cap\
\mathbb{B}(\mathcal{T})^{(s\cdot\underline{f})}\right)\ -\ \,
\mathfrak{F}\!\left(\mathfrak{I}(\upsilon_e)\ \cap\
\mathbb{B}(\mathcal{T})^{(s\cdot\underline{f}+1)}\right)\ \biggr)\
\Biggr)\ .
\end{multline*}
\end{proposition}

{\small Recall that one
structure refining the description of a layer
$\mathbb{B}(n)^{(d)}$ of the Boolean lattice $\mathbb{B}(n)$ of
subsets of an $n$-set, $1\leq d\leq
\left\lfloor\tfrac{n}{2}\right\rfloor$, is that of the Johnson
association scheme, see, e.g.,~\cite[\S{}3.2]{BI},\cite[\S{}2.7,
\S{}9.1]{BCN},\cite[\S{}4.2]{D},\cite[\S{}21.6]{MS}. The {\em
Johnson scheme\/} $\mathbf{J}(n,d)$ is the pair
$(\pmb{X},\boldsymbol{\mathcal{R}})$, where
$\pmb{X}:=\mathbb{B}(n)^{(d)}$ with $|\pmb{X}|=\tbinom{n}{d}$, and
$\boldsymbol{\mathcal{R}}:=(\pmb{R}_0,\pmb{R}_1,\ldots,\pmb{R}_d)$
is a partition of $\pmb{X}\times\pmb{X}$, defined by
\begin{equation*}
\pmb{R}_i:=\bigl\{(x,y):\ \partial(x,y):=d-\rho(x\wedge
y)=i\bigr\}\ ,\ \ \ 0\leq i\leq d\ .
\end{equation*}
For any $x,y\in\pmb{X}$ with $\partial(x,y)=k$, and for any
integers $i$ and $j$, $0\leq i,j\leq d$, the~{\em intersection
numbers\/}
\begin{multline*}
\mathtt{p}^{k}_{ij}\ \ :=\ \  \bigl|\{z\in\pmb{X}:\
\partial(z,x)=i,\ \partial(z,y)=j\}\bigr|\\ =\ \
\Bigl|\,\mathbb{B}(n)^{(d)}\ \ \cap\ \ \Bigl(\
\mathfrak{F}\bigl(\mathfrak{I}(x)\cap\mathbb{B}(n)^{(d-i)}\bigr)-
\mathfrak{F}\bigl(\mathfrak{I}(x)\cap\mathbb{B}(n)^{(d-i+1)}\bigr)\
\Bigr)\\ \cap\ \ \Bigl(\
\mathfrak{F}\bigl(\mathfrak{I}(y)\cap\mathbb{B}(n)^{(d-j)}\bigr)-
\mathfrak{F}\bigl(\mathfrak{I}(y)\cap\mathbb{B}(n)^{(d-j+1)}\bigr)\
\Bigr)\Bigr|
\end{multline*}
are the same. We have
\begin{equation*}
\mathtt{p}_{ij}^k=\sum_c\binom{d-k}{c}\binom{k}{d-i-c}\binom{k}{d-j-c}
\binom{n-d-k}{i+j-d+c}\ ;
\end{equation*}
the quantity $\mathtt{n}_i:=\mathtt{p}^0_{ii}= |\{z\in\pmb{X}:\
\partial(z,x)=i\}|$, for any $x\in\pmb{X}$, called the {\em valency\/} of $\pmb{R}_i$, is
$\mathtt{n}_i=\tbinom{d}{i}\tbinom{n-d}{i}$, see,
e.g.,~\cite{Moon}.
}

The family $\mathbf{K}^{\ast}(\mathcal{M})$ of all tope committees
for a simple oriented matroid $\mathcal{M}:=(E_t,\mathcal{T})$ can
be considered as the poset
$\mathbf{I}_{\frac{1}{2}}\bigl(\mathbb{B}(\mathcal{T}),\Upsilon\bigr)$
of relatively \mbox{$\tfrac{1}{2}$-blocking} elements for the {\em
subset\/} $\Upsilon$ of the Johnson scheme
$\mathbf{J}(|\mathcal{T}|,\tfrac{|\mathcal{T}|}{2}):=(\pmb{X},\boldsymbol{\mathcal{R}})$
on the set $\pmb{X}:=\mathbb{B}(\mathcal{T})^{(|\mathcal{T}|/2)}$,
with the partition
$\boldsymbol{\mathcal{R}}:=(\pmb{R}_0,\pmb{R}_1,\ldots,$
$\pmb{R}_{\frac{|\mathcal{T}|}{2}})$ of $\pmb{X}\times\pmb{X}$,
defined by $\pmb{R}_i:=\bigl\{(x,y):\
\partial(x,y):=\tfrac{|\mathcal{T}|}{2}-\rho(x\wedge y)=i\bigr\}$,
for all $0\leq i\leq \tfrac{|\mathcal{T}|}{2}$.

\section{Layers of Tope Committees Containing no Pairs of Opposites}
\label{Section-No-Opposites}

For any element $e\in E_t$, the corresponding positive halfspace
$\mathcal{T}_e^+$ of a~simple oriented matroid
$\mathcal{M}:=(E_t,\mathcal{T})$ contains no pairs of opposites.
On the~other hand, the tope committees containing pairs of
opposites have no applied value. We now describe the structure of
the family
\begin{equation*}
\overset{\circ}{\mathbf{K}}{}^{\ast}(\mathcal{M}):=\{\mathcal{K}^{\ast}\in
\mathbf{K}^{\ast}(\mathcal{M}):\
T\in\mathcal{K}^{\ast}\Longrightarrow-T\not\in\mathcal{K}^{\ast}\}
\end{equation*}
of all the committees for $\mathcal{M}$ which include no pairs of
opposites.

Denote by $\mathbf{O}'(\mathcal{T})$ a graded
meet-sub-semilattice, of rank $\tfrac{|\mathcal{T}|}{2}$, of the lattice~$\mathbb{B}(\mathcal{T})$, defined in the following way: the
elements of $\mathbf{O}'(\mathcal{T})$ are the subsets of topes
without pairs of opposites, which are ordered by inclusion. This
poset is isomorphic to the face poset of the boundary of a
$\tfrac{|\mathcal{T}|}{2}$-dimensional {\em crosspolytope}.

Let $\gcd(\cdot,\cdot)$ denote the greatest common divisor of a pair of integers, and define a Farey
subsequence $\mathcal{F}\bigl(\mathbf{O}^{\prime}(\mathcal{T}),\rho(a)\bigr)$, associated with an element
$a\in\mathbf{O}^{\prime}(\mathcal{T})$, by {\small
\begin{equation*}
\mathcal{F}\bigl(\mathbf{O}^{\prime}(\mathcal{T}),\rho(a)\bigr):=\left(
\frac{\rho(b\wedge a)}{\gcd(\rho(b\wedge a),\rho(b))} \biggm/
\frac{\rho(b)}{\gcd(\rho(b\wedge a),\rho(b))}:\
b\in\mathbf{O}^{\prime}(\mathcal{T})-\{\hat{0}\}\right)\ ,
\end{equation*}
}

If $\rho(a)<\tfrac{|\mathcal{T}|}{2}$ then the sequence
$\mathcal{F}\bigl(\mathbf{O}^{\prime}(\mathcal{T}),\rho(a)\bigr)=
\bigl(\tfrac{h}{k}\in\mathcal{F}_{|\mathcal{T}|/2}:\ \mbox{$h\leq \rho(a)$}\bigr)$, is a member of the family
of the Farey subsequences of the form $\bigl(\tfrac{h}{k}\in\mathcal{F}_n:\ h\leq m\bigr)$, considered, e.g.,
in~\cite{AZ},\cite[\S{}7]{AM}. In the case where $\rho(a)=\tfrac{|\mathcal{T}|}{2}$, relevant to our context,
the sequence $\mathcal{F}\bigl(\mathbf{O}^{\prime}(\mathcal{T}),\rho(a)\bigr)$ is nothing else than
$\mathcal{F}_{|\mathcal{T}|/2}$, the standard Farey sequence of order $\tfrac{|\mathcal{T}|}{2}$.

The family $\overset{\circ}{\mathbf{K}}{}^{\ast}(\mathcal{M})$ can
be viewed as the subposet
\begin{equation*}
\mathbf{I}_{\frac{1}{2}}\bigl(\mathbf{O}^{\prime}(\mathcal{T}),
\Upsilon\bigr):=\left\{b\in\mathbf{O}^{\prime}(\mathcal{T}):\
\rho(b)>0,\
\tfrac{\rho(b\wedge\upsilon_e)}{\rho(b)}>\tfrac{1}{2},\ \forall
e\in E_t \right\}
\end{equation*}
of all relatively $\tfrac{1}{2}$-blocking elements for the
antichain $\Upsilon$ in $\mathbf{O}^{\prime}(\mathcal{T})$; for
any $k$, $1\leq k\leq\tfrac{|\mathcal{T}|}{2}$, define the
antichain
$\mathbf{I}_{\frac{1}{2},k}\bigl(\mathbf{O}^{\prime}(\mathcal{T}),
\Upsilon\bigr):=\mathbf{O}^{\prime}(\mathcal{T})^{(k)}\cap
\mathbf{I}_{\frac{1}{2}}\bigl(\mathbf{O}^{\prime}(\mathcal{T}),
\Upsilon\bigr)$. This subposet can be described, in view
of~\cite[Theorem~8.4]{AM}, in the following way:


\begin{theorem}
\label{prop:9}
Let $\mathcal{M}$ be a simple oriented matroid, which is not
acyclic, on the ground set $E_t$, with set of topes $\mathcal{T}$.

On the one hand,
\begin{equation*}
\mathbf{I}_{\frac{1}{2}}\bigl(\mathbf{O}^{\prime}(\mathcal{T}),
\Upsilon\bigr) = \dot{\bigcup_{3\leq k\leq\frac{|\mathcal{T}|}{2}}}
\biggl(\
\mathbf{O}^{\prime}(\mathcal{T})^{(k)}\ \, \cap\ \bigcap_{e\in E_t
}\mathfrak{F}\!\left(\mathfrak{I}(\upsilon_e)\ \cap\
\mathbf{O}^{\prime}(\mathcal{T})^{(\lceil(k+1)/2\rceil)}\right)\
\biggr)
\end{equation*}
and, in particular,
\begin{equation*}
\mathbf{I}_{\frac{1}{2},3}\bigl(\mathbf{O}^{\prime}(\mathcal{T}),
\Upsilon\bigr) =
\mathbf{O}^{\prime}(\mathcal{T})^{(3)}\ \, \cap\ \bigcap_{e\in E_t
}\mathfrak{F}\!\left(\mathfrak{I}(\upsilon_e)\ \cap\
\mathbf{O}^{\prime}(\mathcal{T})^{(2)}\right)\ .
\end{equation*}
On the other hand,
\begin{multline*}
\mathbf{I}_{\frac{1}{2}}\bigl(\mathbf{O}^{\prime}(\mathcal{T}),
\Upsilon\bigr)\ \ =\ \ \bigcap_{e\in E_t}\ \ \
\bigcup_{f\in\mathcal{F}_{|\mathcal{T}|/ 2}:\
\frac{1}{2}<f}\ \ \
\bigcup_{s\in\left[1,\left\lfloor|\mathcal{T}| /
(2\overline{f})\right\rfloor\right]}
\\ \Biggl(\ \mathbf{O}^{\prime}(\mathcal{T})^{(s\cdot\overline{f})}\ \ \cap
\ \ \biggl(\ \mathfrak{F}\!\left(\mathfrak{I}(\upsilon_e)\ \cap\
\mathbf{O}^{\prime}(\mathcal{T})^{(s\cdot\underline{f})}\right)\
-\ \, \mathfrak{F}\!\left(\mathfrak{I}(\upsilon_e)\ \cap\
\mathbf{O}^{\prime}(\mathcal{T})^{(s\cdot\underline{f}+1)}\right)\
\biggr)\ \Biggr)\ .
\end{multline*}
\end{theorem}

The\hfill three-tope\hfill committees\hfill for\hfill $\mathcal{M}$\hfill composing\hfill the\hfill antichains\\ $\mathbf{I}_{\frac{1}{2},3}\bigl(\mathbb{B}(\mathcal{T}),
\Upsilon\bigr)$ and $\mathbf{I}_{\frac{1}{2},3}\bigl(\mathbf{O}^{\prime}(\mathcal{T}),
\Upsilon\bigr)$, described in Proposition~\ref{prop:8} and Theorem~\ref{prop:9}, are treated and enumerated in~\cite{M-Three}.

{\small For an integer $m>1$, let $\pm[1,m]$ denote the $2m$-set
$\{-m,-(m-1),\ldots,-2,-1,1,$ $2,\ldots,$ $(m-1),m\}$. Fix some
$d\in[1,m]$, and denote by $\pmb{X}$ the family of all the~\mbox{$d$-subsets} $V\subset\pm[1,m]$ such that
\begin{equation}
\label{eq:11} s\in V \ \ \ \Longrightarrow\ \ \ -s\not\in V\ .
\end{equation}

Let $\mathbf{O}(m)$ denote the family of all the subsets
$V\subset\pm[1,m]$, satisfying~(\ref{eq:11}), ordered by
inclusion, and augmented by a greatest element $\hat{1}$. The
least element of $\mathbf{O}(m)$ is the empty subset of
$\pm[1,m]$. See, e.g.,~\cite{BHP,Readdy} on the structural and
combinatorial properties of the {\em face lattice\/}
$\mathbf{O}(m)$ of an $m$-dimensional {\em crosspolytope\/} ({\em
hyperoctahedron}, {\em orthoplex}), dual to that of an
$m$-dimensional {\em hypercube}. The~graded lattice
$\mathbf{O}(m)$, of rank $m+1$, is {\em Eulerian}. The fundamental
properties of Eulerian posets are described
in~\cite[Chapter~3]{St2}. The poset
$\mathbf{O}^{\prime}(m):=\mathbf{O}(m)-\{\hat{1}\}$ is {\em
simplicial}. For $d\in[1,m]$, the $d$th {\em Whitney number of the
second kind\/} $\mathrm{W}_d(\mathbf{O}(m))$, that is the
cardinality of the layer $\mathbf{O}(m)^{(d)}$, is
$\tbinom{m}{d}2^d$.

Let $\pmb{X}:=\mathbf{O}(m)^{(d)}$. Define a partition
$\boldsymbol{\mathcal{R}}:=(\pmb{R}_0,\pmb{R}_1,\ldots,\pmb{R}_d)$
of $\pmb{X}\times\pmb{X}$ in the~following way:
\begin{equation*}
\pmb{R}_i:=\bigl\{(x,y):\
\partial(x,y):=d-\rho(x\wedge y)=i\bigr\}\ ,\ \ \ 0\leq i\leq d\ ;
\end{equation*}
here, $\rho$ denotes the rank function on $\mathbf{O}(m)$.

For all elements $x\in\pmb{X}$, the quantities
$\mathtt{n}_i:=|\{z\in\pmb{X}:\ \partial(z,x)=i\}|$ are the~same:
$\mathtt{n}_i=\tbinom{d}{i}\sum_{c=0}^{i}\tbinom{i}{c}
\tbinom{m-d}{c}2^c=\tbinom{d}{i}\sum_{c=0}^{i}\tbinom{i}{c}
\tbinom{m-d+c}{i}$, for any $i$, $0\leq i\leq d$.

Recall that in the case where $d:=m$ and $\mathtt{n}_i=\tbinom{m}{i}$, for any $i$, $0\leq i\leq m$, the~pair
$(\pmb{X},\boldsymbol{\mathcal{R}})$ is an association scheme. Indeed, for any $x,y\in\pmb{X}$ with
$\partial(x,y):=m-\rho(x\wedge y)=k$, and for any integers $i$ and $j$, $0\leq i,j\leq m$, the quantities
\begin{multline*}
\mathtt{p}^{k}_{ij}\ \ :=\ \ \bigl|\{z\in\pmb{X}:\
\partial(z,x)=i,\ \partial(z,y)=j\}\bigr|\\ =\ \
\Bigl|\,\mathbf{O}(m)^{(m)}\ \ \cap\ \ \Bigl(\
\mathfrak{F}\bigl(\mathfrak{I}(x)\cap\mathbf{O}(m)^{(m-i)}\bigr)-
\mathfrak{F}\bigl(\mathfrak{I}(x)\cap\mathbf{O}(m)^{(m-i+1)}\bigr)\
\Bigr)\\ \cap\ \ \Bigl(\
\mathfrak{F}\bigl(\mathfrak{I}(y)\cap\mathbf{O}(m)^{(m-j)}\bigr)-
\mathfrak{F}\bigl(\mathfrak{I}(y)\cap\mathbf{O}(m)^{(m-j+1)}\bigr)\
\Bigr)\Bigr|
\end{multline*}
are the same; the intersection numbers
\begin{equation*}
\begin{split}
\mathtt{p}_{ij}^k&=\sum_c\binom{m-k}{c}\binom{k}{m-i-c}
\binom{i+k-m+c}{m-j-c}\binom{m-k-c}{i+j-m+c}\\&=
\begin{cases}\dbinom{m-k}{\frac{i+j-k}{2}}\dbinom{k}{\frac{i-j+k}{2}},&\text{if $i+j+k$ even}\ ,\\
0,&\text{if $i+j+k$ odd}\end{cases}
\end{split}
\end{equation*}
are those of the {\em Hamming\/} association {\em scheme\/}
$\mathbf{H}(m,2)$, see, e.g.,~\cite[\S{}3.2]{BI},\cite[\S{}2.5,
\S{}9.2]{BCN},\cite[\S{}4.1]{D},\cite{Egawa,Enomoto},\cite[\S{}21.3]{MS}.
The {\em $m$-cube\/}
$\mathbf{H}(m,2):=(\pmb{X},\boldsymbol{\mathcal{R}})$ is the
family~$\pmb{X}$, of cardinality $2^m$, of all words
$\pmb{x}\in\{-1,1\}^m$, together with the partition
$\boldsymbol{\mathcal{R}}:=(\pmb{R}_0,\pmb{R}_1,\ldots,\pmb{R}_m)$
of $\pmb{X}\times\pmb{X}$, defined by
$\pmb{R}_i:=\bigl\{(\pmb{x},\pmb{y}):\
\partial(\pmb{x},\pmb{y}):=|\{l\in[1,m]:\ x_l\neq y_l\}|=i\bigr\}$,
$0\leq i\leq m$.
}

Let $\mathcal{M}:=(E_t,\mathcal{T})$ be a simple oriented matroid.
The family $\overset{\circ}{\mathbf{K}}{}^{\ast}(\mathcal{M})$ of
its tope committees, containing no pairs of opposites, can be
considered as the~poset
$\mathbf{I}_{\frac{1}{2}}\bigl(\mathbf{O}^{\prime}(\mathcal{T}),\Upsilon\bigr)$
of relatively $\tfrac{1}{2}$-blocking elements for the {\em
subset\/} $\Upsilon$ of the~association scheme
$(\pmb{X},\boldsymbol{\mathcal{R}})$ on the set
$\pmb{X}:=\mathbf{O}^{\prime}(\mathcal{T})^{(|\mathcal{T}|/2)}$,
with the~partition
$\boldsymbol{\mathcal{R}}:=(\pmb{R}_0,\pmb{R}_1,\ldots,
\pmb{R}_{\frac{|\mathcal{T}|}{2}})$ of $\pmb{X}\times\pmb{X}$,
defined by $\pmb{R}_i:=\bigl\{(x,y):\
\partial(x,y):=\tfrac{|\mathcal{T}|}{2}-\rho(x\wedge y)=i\bigr\}$,
for all $0\leq i\leq \tfrac{|\mathcal{T}|}{2}$; the parameters of
$(\pmb{X},\boldsymbol{\mathcal{R}})$ are those of the
$\tfrac{|\mathcal{T}|}{2}$-cube
$\mathbf{H}(\tfrac{|\mathcal{T}|}{2},2)$.

\end{document}